\documentclass[11pt]{article}

\usepackage{amssymb}
\usepackage[english]{babel}
\usepackage{graphicx}
\usepackage{latexsym}

\usepackage{amsmath,amsfonts,amssymb,latexsym,graphics,epsfig,url}
\usepackage{color}
\usepackage{amsthm}
\usepackage[english]{babel}
\usepackage{graphicx}
\usepackage{soul}
\usepackage[utf8]{inputenc}

\newtheorem{theorem}{Theorem}[section]

\newtheorem{example}[theorem]{Example}

\newtheorem{remark}[theorem]{Remark}

\DeclareMathOperator{\spec}{sp}

\def\Z{\ns Z}

\def\f{\mbox{\boldmath $f$}}

\def\o{\mbox{\rm o}}

\def\vecv{\mbox{\boldmath $v$}}

\def\vec0{\mbox{\boldmath $0$}}

\def\A{\mbox{\boldmath $A$}}
\def\B{\mbox{\boldmath $B$}}

\def\G{\Gamma}

\def\Z{\ns{Z}}

\def\1{\mbox{\boldmath $1$}}

\def\G{\Gamma}

\def\Z{\mathbb Z}

\begin{document}
	\title{On factored lifts of graphs and their spectra
		\thanks{This research has been supported by
			AGAUR from the Catalan Government under project 2021SGR00434 and MICINN from the Spanish Government under project PID2020-115442RB-I00.
			The research of M. A. Fiol was also supported by a grant from the  Universitat Polit\`ecnica de Catalunya with references AGRUPS-2022 and AGRUPS-2023.
			The last two authors acknowledge support from APVV Research Grants 19-0308 and 22-0005, and VEGA Research Grants 1/0567/22 and 1/0069/23.}
	}
	\author{C. Dalf\'o$^a$, M. A. Fiol$^b$,  S. Pavl\'ikov\'a$^c$, and J. \v{S}ir\'a\v{n}$^d$\\
		\\
		{\small $^a$Departament. de Matem\`atica}\\ {\small  Universitat de Lleida, Igualada (Barcelona), Catalonia}\\
		{\small {\tt cristina.dalfo@udl.cat}}\\
		{\small $^{b}$Departament de Matem\`atiques}\\ 
		{\small Universitat Polit\`ecnica de Catalunya, Barcelona, Catalonia} \\
		{\small Barcelona Graduate School of Mathematics} \\
		{\small  Institut de Matem\`atiques de la UPC-BarcelonaTech (IMTech)}\\
		{\small {\tt miguel.angel.fiol@upc.edu}}\\
		{\small $^{c}$ Inst. of Information Engineering, Automation, and Math, FCFT}, \\
		{\small Slovak Technical University, Bratislava, Slovakia}\\
		{\small {\tt sona.pavlikova@stuba.sk} }\\
		{\small $^{d}$Department of Mathematics and Statistics}\\
		{\small The Open University, Milton Keynes, UK}\\
		{\small Department of Mathematics and Descriptive Geometry}\\
		{\small Slovak University of Technology, Bratislava, Slovak Republic}\\
		{\small {\tt  jozef.siran@stuba.sk}}
	}

	\maketitle
	
	\begin{abstract}
		In this note, we introduce the concept of factored lift, associated with a combined voltage graph, as a generalization of the lift graph. We present a new method for computing the eigenvalues and eigenspaces of factored lifts.
	\end{abstract}
	
	\section{Introduction}
	For an informal description of the lifting construction, think of a graph $G$
	(a {\em base graph}) endowed by an assignment $\alpha$ of elements of a group
	$\Gamma$ on arcs of $G$ (a {\em voltage assignment}). The pair $(G,\alpha)$
	gives rise to a {\em lift} $G^\alpha$, a larger graph that may be thought of
	as obtained by `compounding' $|\Gamma|$ copies of $G$, joined in-between in a
	way dictated by the voltage assignment. A general necessary and sufficient
	condition for a graph $H$ to arise as a lift of a (smaller) graph is the existence
	of a subgroup of the automorphism group of $H$ with a free action on vertices
	of $H$; see Gross and Tucker \cite{gt87}. Lifts have found numerous applications
	in areas of graph theory that are as versatile as the degree/diameter problem on the
	one hand and the Map Color Theorem on the other hand.
	In the first case, the diameter of the lift can be conveniently expressed
	in terms of voltages on the edges of the base graph. Besides its theoretical importance,
	this fact can be used to design efficient diameter-checking algorithms, see Baskoro, Brankovi\'c, Miller, Plesn\'{\i}k,  Ryan, and  \v{S}ir\'a\v{n} \cite{bbmprs97}.
	In the context of the degree-diameter problem, we address the interested reader to the comprehensive survey by Miller and  \v{S}ir\'a\v{n} \cite{ms13}.
	Other prominent examples of
	applications also include the use of lift graphs in the study of the automorphisms of $G^{\alpha}$ and in Cayley graphs, which are lifts of one-vertex graphs
	(with loops and semi-edges attached).
	
	Some of the authors, together with Miller and Ryan, introduced a method to obtain the spectrum and eigenvectors of lift graphs (or voltage graphs) in \cite{dfmrs19}. Later, some of us extended this method to digraphs \cite{dfs19}. The following step was to generalize this method to arbitrary lifts of graphs in \cite{dfps21}. Finally, we expanded this generalization to the universal adjacency matrix of arbitrary lifts of graphs in \cite{dfps23}.
	
	In this note, we present the method for obtaining the spectrum and eigenvectors to the more general case when we have voltages in the edges (or arcs) as usual, but now we also have coset voltages associated with the vertices.
	
	The concept of `overlift' was introduced in Reyes, Dalf\'o, Fiol, Messegu\'e \cite{rdfm23}. There, we obtained the spectrum of a lift graph without enough symmetry by changing the polynomial matrix ad hoc but without formal reasoning. Here, we changed the name of the overlight to the more precise `factored lift,' and we gave a complete explanation and its corresponding proof. We think this is a very useful concept 
	because it is a generalization of permutation voltage lifts. There exists another generalization of voltage graphs by Potočnika and Toledo \cite{PoTo21}, where they assign coset voltages to arcs and vertices. In contrast, we assign them to vertices and standard voltages to arcs.
	\smallskip
	
	
	This note is structured as follows. In the rest of this section, we present the definitions related to factored lift graphs and give a pair of examples. In the following section, there is the main result, that is, the theorem that allows us to compute the spectrum of a factored lift graph.
	
	\section{Preliminaries}
	
	Let $\Gamma$ be a graph; as usual, each of its edges is considered to consist of two oppositely directed arcs. (Alternatively, $\Gamma$ can be a directed graph or mixed graph.) We let $V(\Gamma)$ and $A(\Gamma)$ denote the vertex and arc set of $\Gamma$. For a group $G$ and a subgroup $H < G$, we let $G/H$ and $[G:H]$ denote the set of all {\em right} cosets of $H$ in $G$ and the index of $H$ in $G$, respectively.
	\smallskip
	
	A {\em combined voltage assignment} on a graph $\Gamma$ in a group $G$ consists of a pair of functions $(\alpha,\omega)$, where $\alpha$ is an ordinary voltage assignment $\alpha: A(\Gamma)\to G$ (in the case of graphs, with the property that mutually reverse arcs receive mutually inverse voltages), and $\omega$ assigns to every vertex $v \in V(\Gamma)$ a subgroup $\omega(v) < G$. The graph $\G$, together with a combined voltage assignment $(\alpha,\omega)$, is called the {\em combined base graph}.
	\smallskip
	
	A {\em factored lift} $\Gamma^{(\alpha,\omega)}$ of $\Gamma$ for a combined voltage assignment $(\alpha,\omega)$ is the graph (or digraph, or mixed graph) defined as follows. The vertex set of the factored lift is the set $V^{(\alpha,\omega)}= \{(v,H)\ |\  v\in V(\Gamma)\ {\rm and}\ H\in G/\omega(v)\}$. Equivalently, for every vertex $v\in V(\Gamma)$, one has $[G:\omega(v)]$ vertices in the factored lift. For a fixed vertex $v\in V(\Gamma)$ the set $\{ (v,H)\ |\ H\in G/\omega(v)\}$ is the {\em fibre above the vertex $v$}.
	\smallskip
	
	To define the arc set of $\Gamma^{(\alpha,\omega)}$, let $a=uv\in A(\Gamma)$ be an arc emanating from a vertex $u$ and terminating at a vertex $v$, carrying a voltage $\alpha(a) \in G$. For each such arc $a=uv$, there is an arc in the factored lift, emanating from a vertex $(u,H)$ for some $H\in G/\omega(u)$ and terminating at a vertex $(v,K)$ for some $K\in G/\omega(v)$ if and only if $H\alpha(a)\cap K\ne \emptyset$. Such an arc in $\Gamma^{(\alpha,\omega)}$ is denoted $(u,H)_a(v,K)$, and the set of arcs of this form in the factored lift is called the {\em fiber above the arc $a$}.
	\smallskip
	
	Note that if $a=uv$ as above, then for every coset $H\in G/\omega(u)$ there is at least one coset $K\in G/\omega(v)$ such that $H\alpha(a)\cap K\ne \emptyset$, but in general there may be more than one of these cosets. A fiber above an arc $a$ in the factored lift may contain more than one arc emanating from the same vertex. For instance, suppose that $G=\Z_{12}$, $H=3\Z_{12}=\{0,3,6,9\}$, $K=4\Z_{12}=\{0,2,4\}$, and $\alpha(uv)=0$. Then, in the factored graph, vertex $(u,H)$ is adjacent to the vertices
	$(v,K)$, $(v,K+3)$, $(v,K+6)$, and $(v,K+9)$.
	
	Also, in the cases of graphs, if $a$ is an arc coming from an undirected edge $e$, that is, if $a=uv$ and $b=vu$ are the two opposite arcs forming $e$, with $\alpha(b)=\alpha(a)^{-1}$, then the arc $(v,K)_b(a,H)$ is opposite to the arc $(u,H)_a(v,K)$ in the factored lift, because the condition $K\alpha(b)\cap H\ne \emptyset$ is equivalent to $H\alpha(a)\cap K\ne \emptyset$ when $\alpha(b)=\alpha(a)^{-1}$.
	\smallskip
	
	Observe that if $\omega$ assigns to every vertex $v\in V(\Gamma)$ a {\em normal} subgroup of $G$, then left multiplication by an arbitrary element $g\in G$ induces an automorphism of the factored lift, preserving fibers above every vertex and every arc (edge). This is because the incidence condition $H\alpha(a)\cap K\ne \emptyset$ for right cosets $H\in G/\omega(u)$ and $K\in G/\omega(v)$ turns after left multiplication by any $g\in G$ into $gH\alpha(a)\cap gK\ne \emptyset$, which is just another incidence condition for right cosets since $gH=Hg$ and $gK=Kg$, by normality of $H$ and $K$ in $G$.
	\medskip
	
	Let us show a pair of examples.
	
\begin{figure}[t]
		\begin{center}
			\includegraphics[width=10cm]{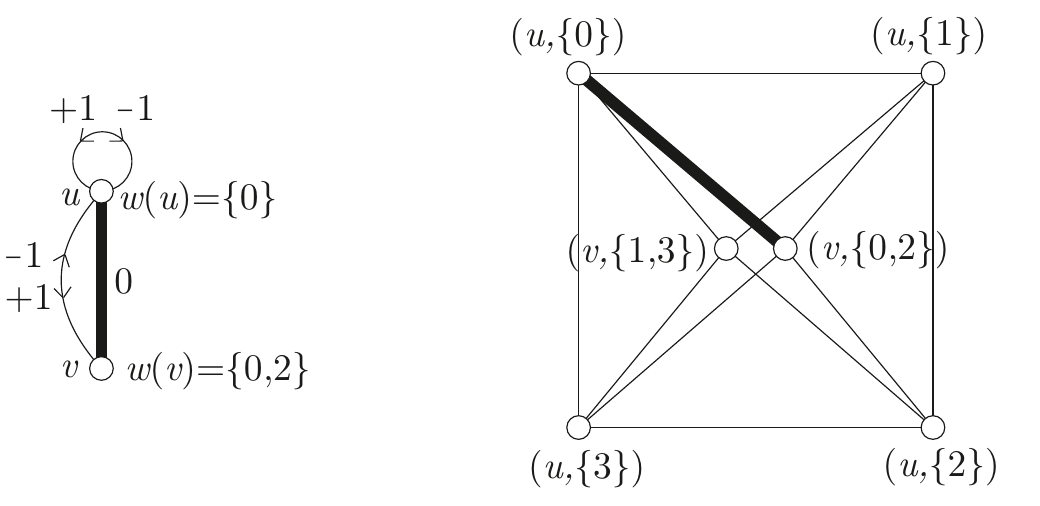}
			\caption{The Johnson graph $J(4,2)$ (or octahedron graph) as a factored lift of the combined base graph on $\Z_4$.}
			\label{fig1}
		\end{center}
	\end{figure}

	\begin{example}
		The representation of the Johnson graph $J(4,2)$ of Figure \ref{fig1}  is an example of a factored lift, with voltage assignment as in the figure; the group is $G=\mathbb{Z}_4=\{0,1,2,3\}$, and if $u$ is the top vertex and $v$ the bottom one, then $\omega(u)$ is the trivial group and $\omega(v)=2\mathbb{Z}_4=\{0,2\}$.    
	\end{example}
	
\begin{figure}[t]
		\begin{center}
			\includegraphics[width=14cm]{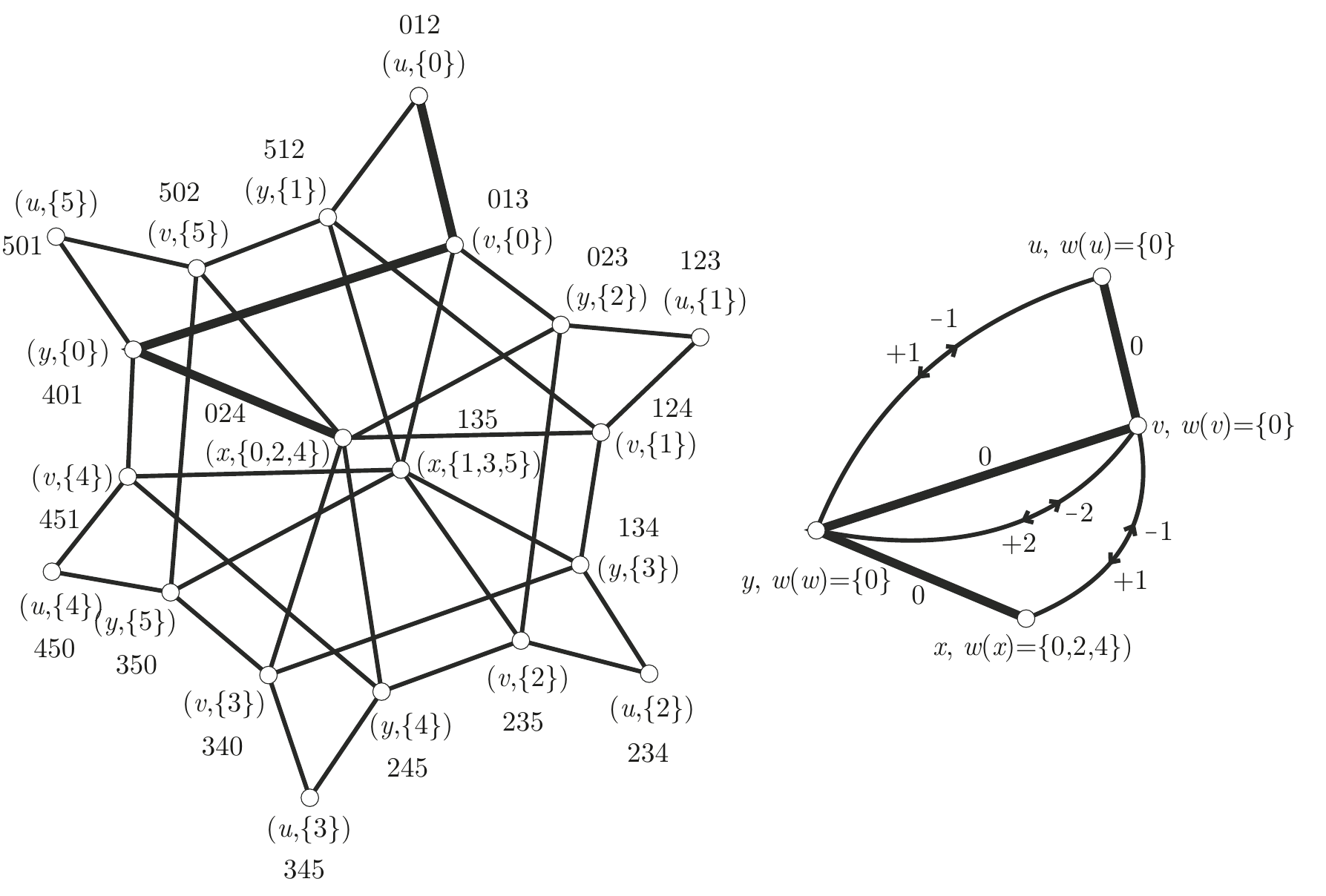}
			\caption{The 3-token graph of $C_6$ and its construction as a factored lift with combined base graph on $\Z_6$.}
			\label{F3(C6)}
		\end{center}
	\end{figure}

	\begin{example}
		\label{ex:F3(C6)}
		Consider now the case of the token graph $F_3(C_6)$ shown in Figure \ref{F3(C6)} (for the definition of token graphs and some of its properties, see Audenaert, Godsil, Royle, and  Rudolph \cite{agrr07}, Fabila-Monroy, Flores-Pe\~{n}aloza, Huemer, Hurtado, Urrutia, and Wood \cite{ffhhuw12}, and Dalf\'o, Duque, Fabila-Monroy, Fiol, Huemer, Trujillo-Negrete, and Zaragoza Mart\'inez \cite{ddffhtz21}).
		Since $3|6$, we have 3 orbits with 6 vertices and one orbit with 2 vertices.
		Then, $F_3(C_6)$ can be obtained as a factored lift with a combined base graph that is a path graph on four vertices, say $u,v,y,x$. Each of these vertices is a representative of one orbit.  For instance, we can take, with the simplified notation $\{i,j,k\}=ijk$ and $(u,v)=uv$, 
		$u=012$, $v=013$, $y=014$, and $x=024$.
		Then, the combined base graph is shown in Figure \ref{F3(C6)}. There, the group is $G=\mathbb{Z}_6$, $\omega(u)=\omega(v)=\omega(y)=\{0\}$, and 
		$\omega(x)=2\mathbb{Z}_6=\{0,2,4\}$. 
		Then, for example, for the coset $H=\{1,3,5\}$
		of $\omega(x)$, there is an edge in the factored lift from the vertex $(x,H)$
		to the vertices $(t,K)$ (where $t\in\{v,y\}$) for the cosets
		$K\in \{ \omega(t)+1,\omega(t)+3,\omega(t)+5 \}$.
	\end{example}
	
	\begin{remark}
		The term `factored lift' is motivated by the fact that it arises by locally factoring the classical (ordinary or
		permutation) lift at vertices (with possibly different groups at different
		vertices). So, a more precise description would be a `locally vertex-factored lift', but we abbreviate this to `factored lift'.
		A factored lift is an ordinary lift if $\omega$ assigns a trivial group (or, more generally, the {\em same normal} subgroup) to every vertex.  
	\end{remark}
	
	\begin{remark}
		Interestingly, if $\omega$ assigns the {\em same but not necessarily normal} subgroup to every vertex, then a factored lift is equivalent to a permutation voltage lift. So, the new concepts look like a potentially useful generalization.
	\end{remark}
	
	
	\section{The spectrum and eigenspaces}
	
	This section shows how to obtain the eigenvalues and eigenspaces of a factored lift from its combined base graph.
	We concentrate on the case when $G$ is a cyclic group $\Z_n$ and, hence, we use its representation with elements $z=\zeta^r=r^{i\frac{r2\pi}{n}}$ for $r=0,1,\ldots, n-1$ (the general case follows similar arguments).
	
	Given a combined base graph $(\Gamma,(\alpha,\omega))$ its {\em associated base graph}
	$(\G,\alpha^+)$ has the same vertices as
	$(\Gamma,(\alpha,\omega))$, all of them associated to the trivial group. 
	In this particular case, it turns out that the factored lift graph $\Gamma^{(\alpha,\omega)}$ is isomorphic to the (standard) lift graph $\G^{\alpha^+}$. The reason is that, as commented in the Introduction, a fiber above an arc $uv$ in the factored lift contains $|\omega(u)|$ 
	arcs emanating from the same vertex.
	
	Let $\B(z)$ be the polynomial matrix of the associated base graph $(\G,\alpha^+)$, where, for each arc $uv$ with voltage $\alpha^+(uv)=i$, the $(u,v)$-entry of $B(z)$ has a term $z^i$.
	
	\begin{theorem}
		\label{th:overlift}
		Let $\Gamma^{(\alpha,\omega)}$ be a factored lift graph associated with the combined voltage graph $\Gamma$ on the vertex set $V(\Gamma)=\{u_1,\ldots,u_{n}\}$ and group $\mathbb{Z}_m$.
		Let $\A$ be the adjacency matrix of  $\Gamma^{(\alpha,\omega)}$, and $\B(z)$ the polynomial $n\times n$ matrix of its associated base graph $(\G,\alpha^+)$. Let 
		$\f= (f_1,\ldots,f_{n})$ be a   $\lambda$-eigenvector of $\B(z)$ with $z=\zeta^r=e^{i\frac{r2\pi}{m}}$ satisfying the following condition:\\
		{\bf P1.}
		Let $o(r)=\frac{m}{\gcd(m,r)}$ and $n_i=[\Z_m:\omega(u_i)]$. For every $u_i\in V(\G)$,  either $f_{i}=0$, or $o(r)$ divides $n_i$.\\
		Then, there exists a corresponding $\lambda$-eigenvector of the factored graph $\Gamma^{(\alpha,\omega)}$.
	\end{theorem}
	
	\begin{proof}
		Let $N=\sum_{i=1}^{n} n_i$
		be the order of the factored lift $\Gamma^{(\alpha,\omega)}$. Let $\f= (f_1,\ldots,f_{n})$ be a $\lambda$-eigenvector of $\B(z)$.
		Then, for every  vertex $X$ of $\Gamma^{(\alpha,\omega)}$, we know that there exist $u_i\in V(\Gamma)$ and $H\in [\Z_m/\omega(u_i)]$, say $H=\omega(u_i)+j$ with $j\in \Z_{n_i}$, such that $X=(u_i,H)$ or, with simplified notation $X=(u_i,j)$. Now, we claim that the vector $\vecv$ with $N$ components of the form
		\begin{equation}
			\label{vector-y}
			v_{X}=v_{(u_i,j)}=f_{i}z^j\qquad u_i\in  V(\Gamma),\  j\in \Z_{n_i} ,
		\end{equation}
		where $z=\zeta^r=e^{ir\frac{2\pi}{m}}$, is a $\lambda$-eigenvector of $\A$, provided that condition {\bf P1} holds.
		To prove it, assume that, in the extended base graph, with the polynomial matrix $\B(z)$, the vertex $u_i$ has $\delta$ out-coming arcs, with voltages $j_1,\ldots,j_{\delta}$,  to the (not necessarily different) vertices $u_{j_1},\ldots,u_{j_{\delta}}$, respectively. Then, the $u_i$-th equality of $\B(z)\f=\lambda\f$ reads
		$$
		\sum_{h=1}^{\delta} f_{j_h} z^{j_h}=\lambda f_i, 
		$$
		which corresponds to the $(u_i,\omega(u_i))$-th equation in $\A\vecv=\lambda\vecv$.
		Now, multiplying all terms by $z^{p}=\zeta^{rp}=e^{ipr\frac{2\pi}{m}}$ 
		with $p=0,1,\ldots, m-1$, we get
		\begin{equation}
			\label{Af=}
			\sum_{h=1}^{\delta} f_{j_h} z^{j_h+p}=\lambda f_i z^{p},\qquad p=0,1,\ldots, m-1,
		\end{equation}
		which should include all equalities (each repeated $n/n_i$ times) in $\A\vecv=\lambda\vecv$ corresponding to the vertices in the same orbit as $u_i$. This holds when, for every $u_i\in V(\G)$, Equation \eqref{Af=} is the same for $p=0$ and $p=n_i$, whence condition {\bf P1} must hold.
		This concludes the proof.
	\end{proof}
	
	\begin{example}
		Here, we show Example \ref{ex:F3(C6)} revisited.
		To find the associated base graph and its voltages, we reason as follows:
		\begin{itemize}
			\item[$(u)$] 
			Vertex $u=012$ is adjacent to $512=y+1$ and $013=v$.\\
			Therefore, $\alpha(uy)=+1$ and $\alpha(uv)=0$.
			\item[$(v)$] 
			Vertex $v=013$ is adjacent to $513=x+1$, $023=y+2$, $012=u$, and $014=y$.\\
			Therefore, $\alpha(vx)=+1$, $\alpha(\underline{vy})=+2$, $\alpha(vu)=0$, and  $\alpha(vy)=0$,\\
			where the underlining represents the second arc between the same pair of vertices. 
			\item[$(y)$] 
			Vertex $y=014$ is adjacent to $514=v-2$, $024=x$, $013=v$, and $015=u-1$. Therefore, $\alpha(\underline{yv})=-2$, $\alpha(yx)=0$,
			$\alpha(yv)=0$, and $\alpha(yu)=-1$.
			\item[$(x)$] 
			Vertex $x=024$ is adjacent to $124=v+1$, $245=y+2$, $014=y$, $034=v+3$, $502=v-1$, and $023=y-2$.\\
			Thus, $\alpha(xv)=z+z^3+z^{-1}$ and $\alpha(xy)=1+z^{2}+z^{-2}$.
		\end{itemize}
		Then, the polynomial matrix is
		$$
		\B(z)  =
		\left(
		\begin{array}{cccc}
			0 & 1 & z & 0 \\
			1 & 0 & 1+z^{2} & z\\
			z^{-1} & 1+z^{-2} & 0 & 1 \\
			0 & z^{-1}+z+z^3 & 1+z^2+z^{-2}  & 0\\
		\end{array}
		\right).
		$$
		Notice that, as expected, such a matrix is obtained from the base matrix of the combined base graph of Figure \ref{F3(C6)}, namely
		$$
		\left(
		\begin{array}{cccc}
			0 & 1 & z & 0 \\
			1 & 0 & 1+z^{2} & z\\
			z^{-1} & 1+z^{-2} & 0 & 1 \\
			0 & z^{-1} & 1  & 0\\
		\end{array}
		\right)
		$$
		by multiplying the last row by $1+z^2+z^{-2}$.
		
		Then, as $F_3(C_6)$ is a factored lift, its eigenvalues can be obtained from $\B(z)$ provided that condition {\bf P1} of Theorem \ref{th:overlift} is fulfilled. In other words, the spectrum of such a matrix contains some `spurious' 
		eigenvalues, not in the spectrum of $F_3(C_6)$. As shown in Table  \ref{taula:C6}, such additional eigenvalues are the four 0's. Following the reasoning in the proof of Theorem \ref{th:overlift}, the reason is that,
		for $r=1,5$, the $\lambda$-eigenvectors of $\B(\zeta^r)$ for $\lambda=0$
		are $(-\frac{1}{2}(1\pm i\sqrt{3}), 0, 0, 1)^{\top}$ and for $r=2,4$
		such $\lambda$-eigenvectors are  $(\frac{1}{2}(1\mp i\sqrt{3}), 0, 0, 1)^{\top}$. Thus, since the last component is $f_4(=f_x)=1(\neq 0)$ and \st{$\pi(x)=2$} $n_x=[\Z_6:\omega(x)]=2$, we have that  $o(1)=o(5)=6\nmid n_x$ and $o(2)=o(4)=3\nmid n_x$. Consequently, none of the above $0$-eigenvectors  yields an eigenvector of $F_3(C_6)$. In contrast, for $r=0,3$ ($z=\pm 1$), the matrices $\B(\pm 1)$ have also an eigenvalue 0 with the corresponding eigenvector $(\mp 1,0,0,1)^{\top}$. Then, although $f_4\neq 0$, both $\o(0)=1$ and $\o(3)=2$ divides $n_x$ and hence such eigenvectors give a $0$-eigenvector of $F_3(C_6)$.
		
		Summarizing, the spectrum of $F_3(C_6)$ is
		$$
		\spec F_3(C_6)=\{4^{[1]},2^{[4]},0^{[2]},-2^{[4]},-4^{[1]}\},
		$$
		indicating that we are dealing with a bipartite graph.
	\end{example}
	
	\begin{table}[!ht]
		\begin{center}
			\begin{tabular}{|c|cccc| }
				\hline
				$\zeta=e^{i\frac{2\pi}{6}}$, $z=\zeta^r$ & $\lambda_{r,1}$  & $\lambda_{r,2}$  & $\lambda_{r,3}$   & $\lambda_{r,4}$  \\
				\hline\hline
				$\spec(\B(\zeta^0))$ &  4 &  0 & $-2$ & $-2$  \\
				\hline
				$\spec(\B(\zeta^1))=\spec(\B(\zeta^5))$ & 2   & $0^*$ & $-1$  &  $-1$ \\
				\hline
				$\spec(\B(\zeta^2))=\spec(\B(\zeta^4))$ &  1 & 1 & $0^*$ &  $-2$ \\
				\hline
				$\spec(\B(\zeta^3))$  & 2 &  2 &  0 &  $-4$ \\
				\hline
			\end{tabular}
		\end{center}
		\caption{All the eigenvalues of the matrices $\B(\zeta^r)$, which yield the eigenvalues of the 3-token graph $F_3(C_6)$ plus four $0$'s (those marked with `*').}
		\label{taula:C6}
	\end{table}
	
	

\end{document}